\documentclass[a4]{article}
 \font\ess=eufb10 scaled\magstep3 

\font\es=eufm10

\def\gP{\mbox{\es {P}}}
\def\J{\mbox{\es {J}}}

\def\lJ{\mbox{\ess {J}}}
\def\CC{\mbox{\es {C}}}

\def\b1{\mbox{\boldmath $1$}}

\def\R{\mbox{\boldmath $R$}}

\def\Z{\mbox{\boldmath $Z$}}

\def\lgP{\mbox{\ess {P}}}

\begin{document}
\title {\bf Orbit types of the compact Lie group $E_7$ in the \\
        complex Freudenthal vector space $\lgP^C$}
\date{}

\author {Takashi Miyasaka and Ichiro Yokota}
\maketitle
\date{}
\section {\bf Introduction}
\hspace*{5pt} Let $\J$ be the exceptional Jordan algebra over $\R$ and $\J^C$ its complexification. Then the simply connected compact exceptional Lie group $F_4$ acts on $\J$ and $F_4$ has three orbit types which are
$$
        F_4/F_4, \quad F_4/Spin(9), \quad F_4/Spin(8). $$ 
\noindent Similarly the simply connected compact exceptional Lie group $E_6$ acts on $\J^C$ and $E_6$ has five orbit types which are
$$
        E_6/E_6, \quad E_6/F_4, \quad E_6/Spin(10), \quad E_6/Spin(9), \quad
        E_6/Spin(8) $$
\noindent ([6]). In this paper, we determine the orbit types of the simply connected compact exceptional Lie group $E_7$ in the complex Freudenthal vector space $\gP^C$. As results, $E_7$ has seven orbit types which are
$$
        E_7/E_7, \quad E_7/E_6, \quad E_7/F_4, 
        \quad E_7/Spin(11),\quad E_7/Spin(10), \quad E_7/Spin(9),
        \quad E_7/Spin(8). $$  
\section {\bf  Preliminaries}
\hspace*{5pt} Let $\CC$ be the division Cayley algebra and $\J = \{ X \in M(3, \CC) \, | \, X^* = X \}$ be the exceptional Jordan algebra with the Jordan multiplication $X \circ Y$, the inner product $(X, Y)$ and the Freudenthal multiplication $X \times Y$ and let $\J^C$ be the complexification of $\J$ with the Hermitian inner product $\langle X, Y\rangle$. (The definitions of $X \circ Y, \, (X, Y), \, X \times Y$ and $\langle X, Y\rangle$ are found in [2]).  Moreover, let $\gP^C = \J^C \oplus \J^C \oplus C \oplus C$ be the Freudenthal $C$-vector space with the Hermitian inner product $\langle P, Q\rangle$. For $P, Q \in \gP^C$, we can define a $C$-linear mapping $P \times Q$ of $\gP^C$. (The definitions of $\langle P, Q\rangle$  and $P \times Q$ are found in [2]). The complex conjugation in the complexified spaces $\CC^C, \J^C$ or $\gP^C$ is denoted by $\tau$. Now, the simply connected compact exceptional Lie groups $F_4, E_6$ and $E_7$ are defined by
\begin{eqnarray*}
           F_4 &=& \{ \alpha \in \mbox{Iso}_R(\J) \, | \,
                    \alpha(X \circ Y) = \alpha X \circ \alpha Y \},\\
          E_6 &=& \{ \alpha \in \mbox{Iso}_C(\J^C) \, | \,
                    \tau\alpha\tau(X \times Y) = \alpha X \times \alpha Y,
                   \langle\alpha X, \alpha Y\rangle = \langle X, Y\rangle \},\\
          E_7 &=& \{ \alpha \in \mbox{Iso}_C(\gP^C) \, | \,
               \alpha(P \times Q)\alpha^{-1} = \alpha P \times \alpha Q,
               \langle\alpha P, \alpha Q\rangle = \langle P, Q\rangle \} \\ 
              &=& \{ \alpha \in \mbox{Iso}_C(\gP^C) \, | \,
               \alpha(P \times Q)\alpha^{-1} = \alpha P \times \alpha Q,
               \alpha(\tau\lambda) = (\tau\lambda)\alpha \}
\end{eqnarray*}
\noindent (where $\lambda$ is the $C$-linear transformation of $\gP^C$ defined by $\lambda(X, Y, \xi, \eta) = (Y, - X, \eta, - \xi))$, respectively. Then we have the natural inclusion $F_4 \subset E_6  \subset E_7$, that is,
\begin{eqnarray*}
       E_6 &=& \{ \alpha \in E_7 \, | \, \alpha(0, 0, 1, 0) = (0, 0, 1, 0) \}
               \subset E_7, \\
       F_4 &=& \{ \alpha \in E_6 \, | \, \alpha E = E \} \subset E_6
               \subset E_7,
\end{eqnarray*}
\noindent where $E$ is the $3 \times 3$ unit matrix. The groups $F_4, E_6$ and $E_7$ have subgroups
\begin{eqnarray*}
Spin(8) &=& \{ \alpha \in F_4 \, | \, \alpha E_k = E_k, k = 1, 2, 3 \} 
                \subset F_4 \subset E_6  \subset E_7,\\
      Spin(9) &=& \{ \alpha \in F_4 \, | \, \alpha E_1 = E_1 \}
               \subset F_4 \subset E_6  \subset E_7, \\
      Spin(10) &=& \{ \alpha \in E_6 \, | \, \alpha E_1 = E_1 \}
               \subset E_6  \subset E_7, \\
      Spin(11) &=& \{ \alpha \in E_7 \, | \, \alpha(E_1, 0, 1, 0) = (E_1, 0, 1,                      0) \} \subset E_7, 
\end{eqnarray*} 
\noindent (the last fact $(E_7)_{(E_1, 0, 1, 0)} = Spin(11)$ will be proved in Theorem 6.(4)), where $E_k$ is the usual notation in $\J^C$, e.g.
      $ E_1 = \pmatrix{1 & 0 & 0 \cr
                      0 & 0 & 0 \cr
                      0 & 0 & 0} $ ([2]).
\section {Orbit types of $F_4$ in $\lJ$ and $E_6$ in $\lJ^C$}
\hspace*{5pt} We shall review the results of orbit types of $F_4$ in $\J$ and $E_6$ in $\J^C$.
\vspace{2mm}

\hspace*{5pt} {\bf Lemma 1} ([1]). {\it Any element $X \in \J$ can be transformed to a diagonal form by some $\alpha \in F_4$} :
\begin{center}
     $ \alpha X = \pmatrix{\xi_1 & 0 & 0 \cr
                          0 & \xi_2 & 0 \cr 
                          0 & 0 & \xi_3}, \quad \xi_k \in \R,    $
      \quad (which is briefly written by $(\xi_1, \xi_2, \xi_3))$.
\end{center}
\noindent {\it The order of $\xi_1, \xi_2, \xi_3$ can be arbitrarily exchanged
under the action of $F_4$.}
\vspace{1mm}

\hspace*{5pt} {\bf Theorem 2} ([6]). {\it The orbit types of the group} $F_4$ {\it in} $\J$ {\it are as follows}.
\vspace{0.5mm}

\hspace*{5pt} (1) {\it The orbit through} $(\xi, \xi, \xi)$ {\it is} $F_4/F_4.$
\vspace{0.5mm}

\hspace*{5pt} (2) {\it The orbit through} $(\xi_1, \xi, \xi)$ ({\it where} $\xi_1 \neq \xi$) {\it is} $F_4/Spin(9).$                                  
\vspace{0.5mm}

\hspace*{5pt} (3) {\it The orbit through} $(\xi_1, \xi_2, \xi_3)$ ({\it where} $\xi_1, \xi_2, \xi_3$ {\it are distinct}) {\it is} $F_4/Spin(8).$  
\vspace{2mm}

\hspace*{5pt} {\bf Lemma 3} ([2],[4]). {\it Any element $X \in \J^C$ can be transformed to the following diagonal form by some $\alpha \in E_6$ }:
\begin{center}
     $ \alpha X = \pmatrix{r_1a & 0 & 0 \cr
                          0 & r_2a & 0 \cr 
                          0 & 0 & r_3a}, \quad r_k \in \R,\,  a \in C, |a| = 1$
\end{center}
\noindent (which is briefly written by $(r_1a, r_2a, r_3a))$.{\it The order of $r_1, r_2, r_3$ can be arbitrarily exchanged
under the action of $E_6$.}
\vspace{2mm}

\hspace*{5pt} {\bf Theorem 4} ([6]). {\it The group} $E_6$ {\it has the following five orbit type in} $\J^C$ :
$$
        E_6/E_6, \quad E_6/F_4, \quad E_6/Spin(10), \quad E_6/Spin(9), \quad
        E_6/Spin(8) $$

\noindent {\it More details,} ($(r_1a, r_2a, r_3a)$ {\it is considered up to a constant} ),
\vspace{0.5mm}

\hspace*{5pt} (1) {\it The orbit through} $(0, 0, 0)$ {\it is} $E_6/E_6.$
\vspace{0.5mm}

\hspace*{5pt} (2) {\it The orbit through} $(1, 1, 1)$ {\it is} $E_6/F_4.$          
\vspace{0.5mm}

\hspace*{5pt} (3) {\it The orbit through} $(1, 0, 0)$ {\it is} $E_6/Spin(10).$                                  
\vspace{0.5mm}

\hspace*{5pt} (4) {\it The orbit through} $(1, r, r)$ ({\it where} $0 < r, r \neq 1)$ {\it is} $E_6/Spin(9).$  
\vspace{0.5mm}

\hspace*{5pt} (5) {\it The orbit through} $(r, s, t)$ ({\it where} $0 \leq r, 0 \leq s, 0 \leq t$ {\it and} $r, s, t$ {\it are distinct}) {\it is} $E_6/Spin(8)$.

\section {Orbit types of $E_7$ in $\lgP^C$}
\hspace*{5pt} {\bf Lemma 5} ([2]). {\it Any element $P \in \gP^C$ can be transformed to the following diagonal form by some $\alpha \in E_7$} :
$$
  \alpha P = (\pmatrix{ar_1 & 0 & 0 \cr
                       0 & ar_2 & 0 \cr
                       0 & 0 & ar_3  },
              \pmatrix{br_1 & 0 & 0 \cr
                       0 & br_2 & 0 \cr
                       0 & 0 & br_3  }, ar, br), \quad
   \begin{array}{l} r_k, r  \in \R, 0 \leq r_k, 0 \leq r, \\
                    a, b \in C, |a|^2 + |b|^2 = 1.
   \end{array}$$
\noindent {\it Moreover, any element $P \in \gP^C$ can be transformed to the following diagonal form by some $\varphi(A)\alpha \in \varphi(SU(2))E_7$} :
\begin{center}
   $ \varphi(A)\alpha P = (\pmatrix{r_1 & 0 & 0 \cr
                       0 & r_2 & 0 \cr
                       0 & 0 & r_3  },
              \pmatrix{0 & 0 & 0 \cr
                       0 & 0 & 0 \cr
                       0 & 0 & 0  }, r, 0), \quad
                     r_k, r  \in \R, 
                     0 \leq r_k, 0 \leq r, $
\end{center}
\noindent (which is briefly written by $(r_1, r_2, r_3; r)$), {\it where $\varphi(A) \in \varphi(SU(2)) \subset E_8$ and commutes with any element $\alpha \in E_7$.}  {\it The order of $r_1, r_2, r_3, r$ can be arbitrarily exchanged
under the action of $E_7$.} (As for the definitions of the groups $E_8$ and $\varphi(SU(2))$, see [2]). The action of $\varphi(A), A \in SU(2),$ on $\gP^C$ is given by
$$
       \varphi(A)P = \varphi\bigl(\pmatrix{a & - \tau b \cr
                                     b & \tau a}\bigr)(X, Y,\xi, \eta) =
   (aX + \tau(bY), aY - \tau(bX), a\xi + \tau(b\eta), a\eta - \tau(b\xi)) .$$  

\hspace*{5pt} {\bf Theorem 6.} {\it The group} $E_7$ {\it has the following seven orbit types in} $\gP^C$ :
$$
        E_7/E_7, \quad E_7/E_6, \quad E_7/F_4, 
        \quad E_7/Spin(11),\quad E_7/Spin(10), \quad E_7/Spin(9),
        \quad E_7/Spin(8). $$  
\noindent {\it More details,}
\vspace{0.5mm}

\hspace*{5pt} (1) {\it The orbit through} (0, 0, 0; 0) {\it is} $E_7/E_7$.
\vspace{0.5mm}

\hspace*{5pt} (2) {\it The orbit through} (0, 0, 0; 1) {\it or} (1, 1, 1; 1) {\it is} $E_7/E_6$.
\vspace{0.5mm}

\hspace*{5pt} (3) {\it The orbit through} (1, 1, 1; 0) {\it or} $(1, 1, 1; r)$ ({\it where} $0 < r, 1 \neq r)$ {\it is} $E_7/F_4$.
\vspace{0.5mm}

\hspace*{5pt} (4) {\it The orbit through} (1, 0, 0; 1) {\it or} $(1, r, r; 1)$ ({\it where} $0 < r, 1 \neq r)$ {\it is} $E_7/Spin(11)$.
\vspace{0.5mm}

\hspace*{5pt} (5) {\it The orbit through} $(1, 0, 0; r)$ ({\it where} $0 < r, 1 \neq r)$ {\it is} $E_7/Spin(10)$.
\vspace{0.5mm}

\hspace*{5pt} (6) {\it The orbit through} $(1, 1, r; 0)$ {\it or} $(1, 1, r; s)$ ({\it where} $0 < r, 0 < s$ {\it and} $1, r, s$ {\it are distinct} ) {\it is} $E_7/Spin(9)$.

\hspace*{5pt} (7) {\it The orbit through} $(1, r, s; 0)$ {\it or} $(1, r, s; t)$ ({\it where} $r, s, t$ {\it are positive and} $1, r, s, t$ {\it are distinct} ) {\it is} $E_7/Spin(8)$.
\vspace{2mm}

\hspace*{5pt} {\bf Proof.}  From Lemma 5, the representatives of orbit types (up to a constant) can be given by the following. 
$$ 
\begin{array}{llll}
   (0, 0, 0; 0), & \quad  (0, 0, 0; 1), & \quad  (0, 0, 1; 1), & \quad (0, 0, 1; r),  \\
   (0, 1, 1; 1), & \quad  (0, 1, 1; r), &  \quad (0, 1, r; s), & \quad (1, 1, 1; 1),  \\
   (1, 1, 1; r),  & \quad (1, 1, r; r),  & \quad (1, 1, r; s), & \quad (1, r, s; t)
\end{array} $$
\noindent where $r, s, t$ are positive, $0, 1, r, s, t$ are distinct and the order of $0, 1, r, s, t$ can be arbitrarily exchanged.
\vspace{1mm}

\hspace*{5pt} (1) The isotropy subgroup $({E_7})_{(0, 0, 0; 0)}$ is obviously $E_7$. Therefore the orbit through (0, 0, 0; 0) is $E_7/E_7$.
\vspace{1mm}

\hspace*{5pt} (2) The isotropy subgroup $({E_7})_{(0, 0, 0; 1)}$ is $E_6$. Therefore the orbit through $(0, 0, 0; 1)$ is $E_7/E_6$.
\vspace{1mm}

\hspace*{5pt} $(2')$ The isotropy subgroup $({E_7})_{(1, 1, 1; 1)}$ is conjugate to $E_6$ in $E_7$. In fact, we know that the following realization of the homogeneous space $E_7/E_6$ :  $E_7/E_6 = \{ P \in \gP^C \, | \, P \times P = 0, \langle P, P \rangle = 1 \} = \mbox{\es M}$ ([4]). Since $\frac{1}{2\sqrt{2}}(E, E, 1, 1)$ and $(0, 0, 1, 0) \in \mbox{\es M}$, there exists $\delta \in E_7$ such that
$$
          \delta\bigl(\frac{1}{2\sqrt{2}}(E, E, 1, 1)\bigr) = (0, 0, 1, 0). $$
\noindent Hence the isotropy subgroup $(E_7)_{(E, E, 1, 1)}$ is conjugate to the isotropy subgroup $(E_7)_{(0, 0, 1, 0)}$ in $E_7$ : $(E_7)_{(E, E, 1, 1)} \sim (E_7)_{(0, 0, 1, 0)}$. On the other hand, since
$$
      \varphi(\pmatrix{\frac{1}{\sqrt{2}} & \frac{1}{\sqrt{2}} \cr
                    - \frac{1}{\sqrt{2}} & \frac{1}{\sqrt{2}}})
                    (E, 0, 1, 0) = \frac{1}{\sqrt{2}}(E, E, 1, 1), $$
\noindent we have $({E_7})_{(E, 0, 1, 0)} = ({E_7})_{(E, E, 1, 1)} \sim (E_7)_{(0, 0, 1, 0)} = E_6 $. Therefore the orbit through $(1, 1, 1; 1)$ is $E_7/E_6$. 
\vspace{1mm}

\hspace*{5pt} (3) The isotropy subgroup $({E_7})_{(1, 1, 1; 0)}$ is $F_4$. In fact, for $\alpha \in E_7$ and $P \in \gP^C$, we have $\alpha(\tau\lambda((P \times P)P)) = \tau\lambda(\alpha((P \times P)P)) = \tau\lambda(\alpha(P \times P)\alpha^{-1}\alpha P) = \tau\lambda((\alpha P \times \alpha P)\alpha P)$. Now, let $P = (1, 1, 1; 0)$. Since $\tau\lambda((P \times P)P) = \frac{3}{2}(0, 0, 0; 1)$, if $\alpha \in E_7$ satisfies $\alpha P = P$, then $\alpha$ also satisfies $\alpha(0, 0, 0; 1) = (0, 0, 0; 1)$. Hence $\alpha \in E_6$, so together with $\alpha E = E$, we have $\alpha \in F_4$. Therefore the orbit through $(1, 1, 1; 0)$ is $E_7/F_4$.
\vspace{1mm}

\hspace*{5pt} $(3')$ The isotropy subgroup $({E_7})_{(1, 1, 1; r)}$ is $F_4$. In fact, let $P = (1, 1, 1; r)$. Since $\tau\lambda((P \times P)P) = \frac{3}{2}(r, r, r; 1)$, if $\alpha \in E_7$ satisfies $\alpha P = P \cdots$ (i), then $\alpha$ also satisfies $\alpha(r, r, r; 1) = (r, r, r; 1) \cdots$ (ii). Take (i) $-$ (ii), then we have $\alpha(1 - r, 1 - r, 1 - r; r - 1) = (1 - r, 1 - r, 1 - r; r - 1)$. Since $1 - r \neq 0$, we have $\alpha(1, 1, 1; - 1) = (1, 1, 1; - 1)$. Together with $\alpha P = P$, we have $\alpha(0, 0, 0; 1) = (0, 0, 0; 1)$ and $\alpha(1, 1, 1; 0) = (1, 1, 1; 0)$. Hence $\alpha \in E_6$ and hence $\alpha \in F_4$. Therefore the orbit through $(1, 1, 1; r)$ is $E_7/F_4$.
\vspace{1mm}

\hspace*{5pt} (4) The isotropy subgroup $({E_7})_{(1, 0, 0; 1)}$ is $Spin(11)$. In fact, we know that $(E_7)_{(E_1, E_1, 1, 1)} = Spin(11)$ ([3]). On the other hand, since
$$
      \varphi(\pmatrix{\frac{1}{\sqrt{2}} & \frac{1}{\sqrt{2}} \cr
                    - \frac{1}{\sqrt{2}} & \frac{1}{\sqrt{2}}})
                    (E_1, 0, 1, 0) = \frac{1}{\sqrt{2}}(E_1, E_1, 1, 1), $$
\noindent we have $({E_7})_{(E_1, 0, 1, 0)} = ({E_7})_{(E_1, E_1, 1, 1)} = Spin(11)$. therefore the orbit through $(1, 0, 0; 1)$ is $E_7/Spin(11)$. Therefore the orbit through $(1, 0, 0; 1)$ is $E_7/Spin(11)$.
\vspace{1mm}

\hspace*{5pt} $(4')$ The isotropy subgroup $({E_7})_{(1, r, r; 1)}$ is $Spin(11)$. In fact, let $P = (1, r, r; 1)$. Since $\tau\lambda((P \times P)P) = \frac{3}{2}(r^2, r, r; r^2)$, if $\alpha \in E_7$ satisfies $\alpha P = P \cdots$ (i), then $\alpha$ also satisfies $\alpha(r^2, r, r; r^2) = (r^2, r, r; r^2) \cdots$ (ii). Take (i) $-$ (ii), then we have $\alpha(1 - r^2, 0, 0; 1 - r^2) = (1 - r^2, 0, 0; 1 - r^2)$. Since $1 - r^2 \neq 0$, we have $\alpha(1, 0, 0, 1) = (1, 0, 0; 1)$. Hence $\alpha \in Spin(11)$. Therefore the orbit through $(1, r, r; 1)$ is $E_7/Spin(11)$.
\vspace{1mm}

\hspace*{5pt} (5) The isotropy subgroup $({E_7})_{(1, 0, 0; r)}$ is $Spin(10)$. In fact, for $\alpha \in E_7$ and $P \in \gP^C$, we have $\alpha((P \times P)\tau\lambda P) = (\alpha(P \times P) \alpha^{-1})\alpha(\tau\lambda P) = (\alpha P \times \alpha P)\tau\lambda(\alpha P)$. Now, let $P = (1, 0, 0; r).$  Since $(P \times P)\tau\lambda P = - \frac{1}{2}(r^2, 0, 0; r)$, if $\alpha \in E_7$ satisfies $\alpha P = P \cdots$ (i), then $\alpha$ also satisfies $\alpha(r^2, 0, 0; r) = (r^2, 0, 0; r) \cdots$ (ii). Take (i) $-$ (ii), then we have $\alpha(1 - r^2, 0, 0; 0) = (1 - r^2, 0, 0; 0)$. Since $1 - r^2 \neq 0$, we have $\alpha(1, 0, 0; 0) = (1, 0, 0; 0) \cdots$ (iii). Take (i) $-$ (iii), then $\alpha(0, 0, 0; r) = (0, 0, 0; r)$, that is, $\alpha(0, 0, 0; 1) = (0, 0, 0; 1)$. Hence $\alpha \in E_6$ and $\alpha E_1 = E_1$. Thus $\alpha \in Spin(10)$. Therefore the orbit through $(1, 0, 0; r)$ is $E_7/Spin(10)$.\vspace{1mm}

\hspace*{5pt} (6) The isotropy subgroup $({E_7})_{(1, 1, r; 0)}$ is $Spin(9)$. In fact, let $P = (1, 1, r; 0)$. Since $\tau\lambda((P \times P)P) = \frac{3}{2}(0, 0, 0; r)$, if $\alpha \in E_7$ satisfies $\alpha P = P$, then $\alpha$ also satisfies $\alpha(0, 0, 0; 1) = (0, 0, 0; 1)$. Hence $\alpha \in E_6$, so together with $\alpha P = P$, we have $\alpha \in Spin(9)$ (Theorem 4.(4)). Therefore the orbit through $(1, 1, r; 0)$ is $E_7/Spin(9)$.
\vspace{1mm}

\hspace*{5pt} $(6')$ The isotropy subgroup $({E_7})_{(1, 1, r; s)}$ is $Spin(9)$. In fact, let $P = (1, 1, r; s)$. Since $\tau\lambda((P \times P)P) = \frac{3}{2}(rs, rs, s; r)$, if $\alpha \in E_7$ satisfies $\alpha P = P \cdots$ (i), then $\alpha$ also satisfies $\alpha(rs, rs, s; r) = (rs, rs, s; r) \dots$ (ii). Take (i) $\times$ $r$ $-$ (ii) $\times$ $s$,  then we have $\alpha(r(1 - s^2), r(1 - s^2), r^2 - s^2; 0) = (r(1 - s^2)s, r(1 - s^2), r^2 - s^2; 0)$. Since $r(1 - s^2), r^2 - s^2$ are non-zero and $r(1 - s^2) \neq r^2 - s^2$, from (6) we have $\alpha \in Spin(9)$. Therefore the orbit through $(1, 1, r; s)$ is $E_7/Spin(9)$.
\vspace{1mm}

\hspace*{5pt} (7) The isotropy subgroup $({E_7})_{(1, r, s; 0)}$ is $Spin(8)$. In fact, let $P = (1, r, s; 0)$. Since $\tau\lambda((P \times P)P) = \frac{3}{2}(0, 0, 0; rs)$, if $\alpha \in E_7$ satisfies $\alpha P = P$, then $\alpha$ also satisfies $\alpha(0, 0, 0; 1) = (0, 0, 0; 1)$. Hence $\alpha \in E_6$, so together with $\alpha P = P$, we have $\alpha \in Spin(8)$ (Theorem 4.(5)). Therefore the orbit through $(1, r, s; 0)$ is $E_7/Spin(8).$
\vspace{1mm}

\hspace*{5pt} $(7')$ The isotropy subgroup $({E_7})_{(1, r, s; t)}$ is $Spin(8)$. In fact, let $P = (1, r, s; t)$. Since $\tau\lambda((P \times P)P) = \frac{3}{2}(rst, st, rt; rs)$, if $\alpha \in E_7$ satisfies $\alpha P = P \cdots$ (i), then $\alpha$ also satisfies $\alpha(rst, st, rt; rs) = (rst, st, rt; rs) \cdots$ (ii). Take (i) $\times$ $rs$ $-$ (ii) $\times$  $t$,  then we have $\alpha(rs(1 - t^2), s(r^2 - t^2), r(s^2 - t^2);  0) = (rs(1 - t^2), s(r^2 - t^2), r(s^2 - t^2);  0)$. Since $rs(1 - t^2), s(r^2 - t^2)$ and $r(s^2 - t^2)$ are non-zero and distinct, from (7) we have $\alpha \in Spin(8)$. Therefore the orbit through $(1, r, s; t)$ is $E_7/Spin(8)$.
\vspace{2mm}

\begin{center} {\bf References}
\end{center}
\vspace{2mm}

 [1] H. Freudenthal, {\it Oktaven, Ausnahmegruppen und Oktavengeometrie,} Math. Inst. Rijksuniv. te Utrecht, 1951.
\vspace{1mm}

 [2] T. Miyasaka, O. Yasukura and I. Yokota, {\it Diagonalization of an element $P$ of $\gP^C$ by the compact Lie group $E_7$}, Tsukuba J. Math., to appear.    \vspace{1mm}

 [3] O. Yasukura and I. Yokota, {\it Subgroup $(SU(2) \times Spin(12))/\Z_2$ of compact simple Lie group $E_7$ and non-compact simple Lie group $E_{7, \sigma}$ of type $E_{7(-5)}$,}  Hiroshima Math. J. 12(1982), 59-76.                      \vspace{1mm}

 [4] I. Yokota, {\it Simply connected compact simple Lie group $E_{6(-78)}$ of type $E_6$ and its involutive automorphisms}, J. Math., Kyoto Univ., 20(1980), 447-473.
\vspace{1mm}

 [5]  I. Yokota, {\it  Realization of involutive automorphisms $\sigma$ of exceptional Lie groups $G$,} part II, $G = E_7$, Tsukuba J. Math., 14(1990), 379 - 404.
\vspace{1mm}

 [6] I. Yokota, {\it Orbit types of the compact Lie group $E_6$ in the complex exceptional Jordan algebra $\J^C$}, Inter. Symp. on nonassociative algebras and related topics, Hiroshima, Japan, World Scientific,1990, 353-359.
\vspace{1mm}

 [7] I. Yokota, {\it Exceptional simple Lie groups} (in Japanese), Gendai-Sugakusha, (1992).


\end{document}